\documentclass{article}
\usepackage{college-math-j}
\usepackage{bm}

\theoremstyle{theorem}
\newtheorem{theorem}{Theorem}

\theoremstyle{definition}

\begin{document}

\title{Rebalance your portfolio without selling}
\author{Jay Bartroff} 

\maketitle


\begin{biog} 
\item[\biogpic{\includegraphics[width=84pt]{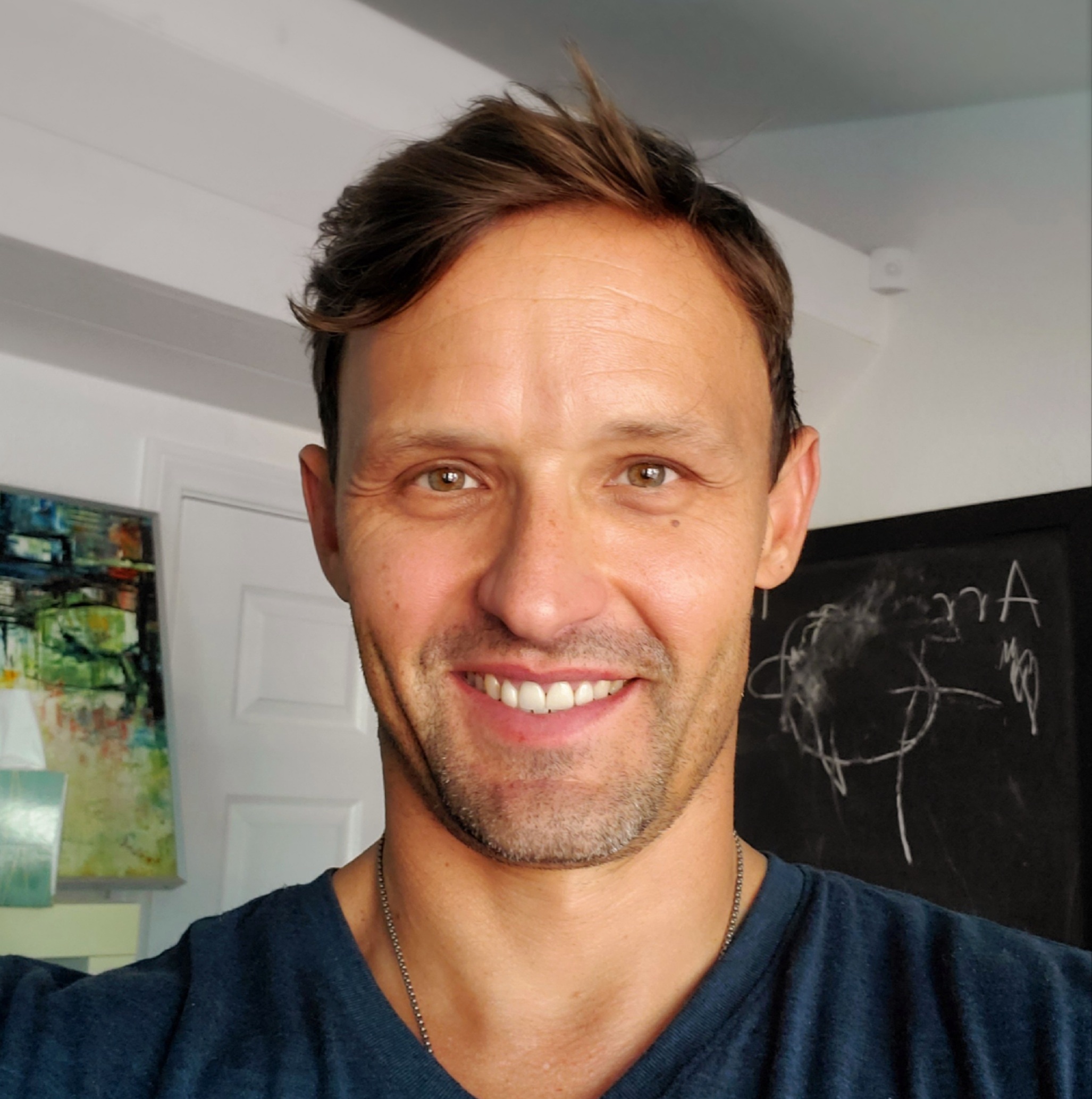}}Jay Bartroff] (bartroff@austin.utexas.edu) received his bachelor's degree from U.C. Berkeley, Ph.D.\ from Caltech, and is now Professor of Statistics \& Data Sciences at the University of Texas at Austin. He is interested in statistics, probability, and occasionally rebalancing his (modest, ahem) portfolio.
\end{biog}

\noindent How do you bring your assets as close as possible to your target allocation by only investing a fixed amount of additional funds, and not selling any assets?  We look at two versions of this problem which have simple, closed form solutions revealed by basic calculus and algebra.

\section{Portfolio rebalancing without selling}
Suppose you own $n$ different types of assets (stocks, bonds, and cash, say, with $n=3$) with values\footnote{Permitting $x_i=0$ allows the possibility of adding a new asset to your portfolio. To keep the story simple we do not allow $x_i<0$, which in finance lingo is called a \textit{short position} and means selling a borrowed asset, or a promise to buy a certain asset at a later date. But the formulas that follow hold under short positions as well as \textit{long positions}~$x_i\ge 0$.}~$x_1,\ldots,x_n\ge 0$. You have in mind target proportions~$p_1,\ldots,p_n\ge 0$ (with $\sum_{i=1}^n p_i=1$) for your assets but the current values do not necessarily satisfy these; that is, $x_i$ is not necessarily equal to $p_i x$ where $x:=\sum_{i=1}^n x_i$ is the current total value of your portfolio.  When it comes time to add an amount $y>0$ to this portfolio (say, through a paycheck deduction to your retirement savings), what is the best way to allocate these additional funds among the assets to bring your portfolio as close as possible to your target allocations?

We provide two answers to this question, driven by two interpretations of the phrase ``as close as possible.'' The first minimizes the sum of squared deviations from the target allocations, which we call the $\ell^2$ problem. The second minimizes the sum of absolute deviations from the target allocations, which we call the $\ell^1$ problem. Although numerical methods abound for solving problems of this type, both happen to have remarkably simple closed-form solutions revealed by only basic calculus and algebra; in one case we do make use of foundational results from convex optimization but they will remain in the background. 

In addition to rebalancing a financial portfolio, this same problem could arise when a city government allocates additional funding $y$ to its $n$ districts based on their relative populations~$p_i$, additional computing  time~$y$ on a supercomputer is allocated among its $n$ existing users with relative priorities~$p_i$, or an additional amount~$y$ of generated energy is allocated to the electrical grid's $n$ nodes according to their relative demand~$p_i$. But the problem we consider arises even outside of ``resource allocation'' problems such as these, and after the $\ell^2$ problem we discuss one such setting where the same calculation is required to perform statistical estimation on a simplex.

Let us go back to the financial portfolio terminology in the first paragraph. The aversion to selling assets may come from, say, the desire to avoid transaction costs, or the difficulty in selling certain financial instruments before a certain date  (e.g., the maturity date of a savings bond). But if we were open to both buying and selling assets,  it is simple algebra to calculate the adjustments
\begin{equation}\label{del.def}
\delta_i:=p_i(x+y)-x_i
\end{equation}
that would bring your portfolio into the proportions $p_i$. That is, to rebalance the old-fashioned way, the $i$th asset currently valued at $x_i$ should be brought to value $x_i+\delta_i$. We call the $\delta_i$ the \textit{naive adjustments}.    Although $\sum_{i=1}^n\delta_i=y>0$ reflecting that the portfolio experiences a net increase of $y$, this old-fashioned rebalancing may require selling (for $\delta_i$ negative) as well as buying (for $\delta_i$ positive). Thus, to bring these values ``as close as possible'' to the target \emph{without} selling, we cannot blindly use the naive adjustments.

One interpretation of our goal is to find adjustments $y_1,\ldots,y_n\ge 0$ subject to $\sum_{i=1}^n y_i=y$ minimizing the sum 
\begin{equation}\label{ss.dev.0}
\sum_{i=1}^n\left(\frac{x_i+y_i}{x+y}-p_i\right)^2
\end{equation}
of squared deviations of the achieved proportions~$(x_i+y_i)/(x+y)$ from the targets~$p_i$. It will make life easier to multiply the objective function~\eqref{ss.dev.0} by the constant $(x+y)^2$, so after a little algebra we see that an equivalent goal is to
\begin{equation}\label{l2.prob}
\mbox{find $y_1,\ldots,y_n\ge 0$ minimizing}\quad \sum_{i=1}^n\left(y_i -\delta_i \right)^2 \quad\mbox{subject to}\quad \sum_{i=1}^n y_i=y.
\end{equation}
Writing the problem in this way shows that our problem is equivalent to finding nonnegative adjustments~$y_i$ totaling $y$ that are ``as close as possible'' to the naive adjustments~$\delta_i$, in the $\ell^2$ sense because \eqref{l2.prob} uses the square of the $\ell^2$ norm to penalize the distance between these quantities.  Other choices of norm, like the $\ell^1$ norm, may of course be reasonable too and so a more general and compact form of the problem is to
\begin{equation}\label{gen.prob}
\mbox{find $\bm{y} \in \mathbb{R}_+^n$ minimizing}\quad \Vert \bm{y} -\bm{\delta} \Vert \quad\mbox{subject to}\quad \bm{1}'\bm{y}=y,
\end{equation} where $\bm{y}$ and $\bm{\delta}$ are the (column) vectorized versions of those variables, $\bm{1}$ is an $n$-long column vector of $1$'s, $\mathbb{R}_+^n=[0,\infty)^n$ is the nonnegative orthant, and $\Vert\cdot\Vert$ is any norm of choice, or any increasing function thereof. The general problem \eqref{gen.prob} with \mbox{$||\cdot||$} any norm is a \textit{convex optimization problem}, the objective function $\Vert \bm{y} -\bm{\delta} \Vert $ and constraints ($y_i$'s nonnegative and sum to $y$) being convex functions. This affords us the comforting fact that any local minimum, should we be so lucky to find one, is also a global minimum. Convex optimization problems also have powerful numerical techniques\footnote{In addition to our $\ell^2$ problem being a \textit{quadratic program}, it is also a special case of \textit{Markowitz portfolio optimization} \cite{Markowitz52}, originated by Nobel laureate Harry~Markowitz, although this particular version does not appear to have been solved before.}  for solving, however we focus on two versions that can be solved directly.  

\section{The $\ell^2$ problem}\label{sec:L2}

The solution of the $\ell^2$ problem~\eqref{l2.prob} is given in Theorem~\ref{thm:l2}. In it, the naive adjustments~$\delta_i$ (given by \eqref{del.def}) are thresholded at a carefully chosen level~$\lambda^*$, and the optimal adjustments~$y_i^*$ end up being the positive part of the excess  over this threshold, $y_i^*=(\delta_i-\lambda^*)^+$ where $z^+$ denotes the positive part~$\max\{z,0\}$. The number of assets that are added to is given by $k^*$ in \eqref{l1.k.lambda}, and the remaining $n-k^*$  assets whose naive adjustments~$\delta_i$ fall below the threshold~$\lambda^*$ remain unchanged. Recall that $x_i$ and $p_i$ denote the initial amount and target proportion, respectively, of the $i$th asset, $x=\sum_{i=1}^n x_i$ is the initial portfolio value, and $y$ denotes the amount added to the portfolio. In what follows let $[n]$ denote $\{1,2, \ldots,n\}$.

The naive adjustments in the portfolio rebalancing problem have a special property -- they sum to $y$, the same constraint asked of the solutions $y_i^*$. In Theorem~\ref{thm:l2} we give the solution to a more general version of the problem where this is not required of the $\delta_i$, but they can be arbitrary real numbers although, to keep the notation simple, we assume they are pre-ordered. The theorem can of course be applied to naive adjustments from a portfolio rebalancing problem, after being put in order.

\begin{theorem}\label{thm:l2} Given arbitrary $y>0$ and $\delta_1\ge \delta_2\ge\ldots\ge \delta_n$,  define 
\begin{equation}\label{l1.k.lambda}
\lambda^*:=\frac{\sum_{i=1}^{k^*}\delta_i-y}{k^*}\quad\mbox{where}\quad 
k^*:=\max\left\{k\in[n]: \sum_{i=1}^k(\delta_i-\delta_k)<y\right\}.
\end{equation}
Then $y_i^*=(\delta_i-\lambda^*)^+$, $i\in[n]$, is the unique solution to the $\ell^2$ problem~\eqref{l2.prob}. In particular, $y_1^*,\ldots, y_{k^*}^*>0$ and $y_{k^*+1}^*=\ldots=y_n^*=0$.
\end{theorem}

In the portfolio rebalancing problem, the last sentence of the theorem says that the $k^*$ assets with the largest $\delta_i$'s are added to, while the remaining $n-k^*$ assets are unchanged.

\begin{proof} As is common with constrained optimization problems, we will consider the \textit{Lagrangian dual} problem which explicitly incorporates the constraint $\sum y_i=y$ into the objective function. That is, first we minimize
$$f(\bm{y},\lambda):= \sum_{i=1}^n (y_i-\delta_i)^2+2\lambda\left( \sum_{i=1}^n y_i-y\right) $$ 
over $\bm{y}\in\mathbb{R}_+^n$, obtaining $\bm{y}^*(\lambda)$, and then we \textit{maximize} $g(\lambda):= f(\bm{y}^*(\lambda),\lambda)$ over $\lambda$, obtaining $\lambda^*$. This maximum $g(\lambda^*)$ is equal to the minimum value\footnote{The precise result we rely on here is \textit{Slater's theorem} which says that \textit{strong duality} holds (i.e., the maximum of the Lagragian dual problem equals the minimum of the original problem) if the \textit{refined Slater conditions} hold, which in this case amount to $y_i^*\ge 0$ for  all $i$ and $y_i^*>0$ for some $i$; see Chapter~5.2.3 of \cite{Boyd04}.} of the sum of squares in the original problem~\eqref{l2.prob}, and since the resulting $\bm{y}^*(\lambda^*)$ will turn out to be unique, it is the unique minimizer.
 
 We have
\begin{equation}\label{dLdy}
\frac{\partial f}{\partial y_i}=2(y_i-\delta_i) +2\lambda
\end{equation} and this of course vanishes when $y_i=\delta_i-\lambda$, however this value may be negative.  Thus, to minimize $f$ we take $y_i^*=\delta_i-\lambda$ for $i$ such that $\delta_i>\lambda$ and, since \eqref{dLdy} is increasing in $y_i$, the remaining $y_i^*$ should be taken as close as possible to  $\delta_i-\lambda$ while remaining nonnegative, i.e., $y_i^*=0$ for the other $i$. We write this compactly as $y_i^*=(\delta_i-\lambda)^+$. 

Let $\delta_0:=\infty$ and define $\kappa(\lambda)$ to be the largest $k\in[n]\cup\{0\}$ such that $\delta_k>\lambda$. 
We have $y_i^*=\delta_i-\lambda$ for $\kappa(\lambda)\ge i\in [n]$ and $y_i^*=0$ for $\kappa(\lambda)< i\in[n]$, so plugging this into $f$ we have 
\begin{multline}
\label{g.lam}
g(\lambda)= \sum_{\kappa(\lambda)\ge i\in [n]}\lambda^2+\sum_{\kappa(\lambda)< i\in[n]}\delta_i^2+2\lambda\left(\sum_{\kappa(\lambda)\ge i\in [n]}(\delta_i-\lambda)-y\right) \\
=-\kappa(\lambda)\lambda^2+2\lambda\left(\sum_{\kappa(\lambda)\ge i\in [n]} \delta_i -y\right)+\sum_{\kappa(\lambda)< i\in[n]}\delta_i^2.
\end{multline}
We claim that $g(\lambda)$ is unimodal with mode at $\lambda=\lambda^*$ given by \eqref{l1.k.lambda}, which will prove the theorem.  To prove this, we will show that $g(\lambda)$ is differentiable with 
$$g'(\lambda) \quad \left\{\begin{array}{c}
 >\\
 =\\
 < \end{array}\right\} \quad\mbox{$0$ for $\lambda$}\quad
 \left\{ \begin{array}{c}
 <\\
 =\\
 > \end{array}\right\}\quad
 \lambda^*.$$
Between consecutive values of $\delta_i$, $\kappa(\lambda)$ is constant so \eqref{g.lam} is differentiable there with derivative
\begin{equation}\label{g'} 
g'(\lambda)=-2\kappa(\lambda)\lambda+ 2 \left(\sum_{\kappa(\lambda)\ge i\in[n]} \delta_i -y\right).
\end{equation}
So the only question of differentiability is at the $\delta_i$. First see that $g(\lambda)$ is continuous at the $\delta_i$.  Suppose $\delta_{j+1}<\delta_j=\delta_{j-1}=\dots=\delta_{k+1}<\delta_{k}$ for some $1\le k<j\le n$ (taking $\delta_{n+1}:=-\infty$ to handle the $j=n$ case). Then $\kappa(\delta_j)=k$ so
\begin{multline*}
\lim_{\lambda \rightarrow \delta_j^+}g(\lambda) =g(\delta_j)=-k\delta_j^2+2\delta_j\left(\sum_{k\ge i\in[n]} \delta_i -y\right)+\sum_{k<i\in[n]}\delta_i^2\\
=-j\delta_j^2+(j-k)\delta_j^2+2\delta_j\left(\sum_{j\ge i\in[n]} \delta_i -(j-k) \delta_j-y\right)+\sum_{j<i\in[n]}\delta_i^2+(j-k)\delta_j^2\\
=-j\delta_j^2+ 2\delta_j\left(\sum_{j\ge i\in[n]} \delta_i -y\right)+\sum_{j<i\in[n]}\delta_i^2  =\lim_{\lambda \rightarrow \delta_j^-}g(\lambda).
\end{multline*}
By a similar argument, $g$ is differentiable there as well:
\begin{align*}
\lim_{\lambda \rightarrow \delta_j^+} g'(\lambda) &=-2k\delta_j +2 \left(\sum_{k\ge i\in[n]} \delta_i -y\right) \\
&=-2j\delta_j +2(j-k)\delta_j +2 \left(\sum_{j\ge i\in[n]} \delta_i -(j-k) \delta_j-y\right) \\
&=-2j\delta_j + 2 \left(\sum_{j\ge i\in[n]} \delta_i -y\right)\\
&=\lim_{\lambda \rightarrow \delta_j^-}g'(\lambda).
\end{align*}

Next we claim that $\kappa(\lambda^*)=k^*$. If it were that $\kappa(\lambda^*) <k^*$ then
$$\delta_{k^*}\le \lambda^*=\frac{\sum_{i\le k^*} \delta_i -y}{k^*}\quad\Leftrightarrow \quad \sum_{i\le k^*} (\delta_i -\delta_{k^*}) \ge y,$$
contradicting the definition~\eqref{l1.k.lambda} of $k^*$. On the other hand, if it were that $\kappa(\lambda^*) >k^*$ then
$$\delta_{k^*+1}> \lambda^*=\frac{\sum_{i\le k^*} \delta_i -y}{k^*}\quad\Leftrightarrow \quad y>\sum_{i\le k^*} \delta_i -k^*\delta_{k^*+1}=\sum_{i\le k^*+1} (\delta_i - \delta_{k^*+1}),$$
again contradicting the definition of $k^*$.

With $\kappa(\lambda^*)=k^*$ established, it is clear from \eqref{g'} that $g'(\lambda^*)=0$. For $\lambda<\lambda^*$ we have $\kappa(\lambda)\ge k^*$ so, using \eqref{g'},
\begin{align}
g'(\lambda)&=-2k^*\lambda-2(\kappa(\lambda)-k^*)\lambda+2 \left(\sum_{k^*\ge i\in[n]} \delta_i+\sum_{k^*<i\le \kappa(\lambda)} \delta_i -y\right)\nonumber\\
&=\left[-2k^*\lambda+2 \left(\sum_{k^*\ge i\in[n]} \delta_i -y\right)\right] +2 \sum_{k^*<i\le \kappa(\lambda)} (\delta_i-\lambda),\label{g'.l<l*}
\end{align}
where the sum $\sum_{k^*<i\le \kappa(\lambda)}=0$ if $k^*=\kappa(\lambda)$. In any case, the last term in \eqref{g'.l<l*} is nonnegative because each summand is, by definition of $\kappa(\lambda)$. The expression in square brackets in \eqref{g'.l<l*} is a decreasing function of $\lambda$, hence is greater than $g'(\lambda^*)=0$, showing that $g'(\lambda)$ is positive.

The arguments for showing that $g'(\lambda)<0$ for $\lambda>\lambda^*$ are similar, so we omit them here.
\end{proof}

\section{The $\ell^1$ problem}\label{sec:L1}

So far we have only considered the sum of squared deviations (or square of the $\ell^2$ norm) to penalize how far our allocations are from their targets, but other norms may be reasonable too. An obvious alternative is to replace the sum of squares in \eqref{ss.dev.0} and \eqref{l2.prob} by the sum of absolute deviations (the $\ell^1$ norm) leading to the problem,
\begin{equation}
\label{l1.prob}
\mbox{find $y_1,\ldots,y_n\ge 0$ minimizing}\quad \sum_{i=1}^n|y_i -\delta_i| \quad\mbox{subject to}\quad  \sum_{i=1}^n y_i=y.
\end{equation}
Like the $\ell^2$ problem in Theorem~\ref{thm:l2}, this can be solved directly.  But unlike the $\ell^2$ problem, the solutions to the $\ell^1$ problem are not unique in general, but can be characterized geometrically using a hyperplane. As with Theorem~\ref{thm:l2}, in Theorem~\ref{thm:l1} we solve a slightly more general version of the problem where the values $\delta_i$ are arbitrary, and are not assumed to sum to $y$ as the naive adjustments do in the portfolio rebalancing problem.

\begin{theorem}\label{thm:l1} Let $y>0$ and $\delta_1, \delta_2, \dots, \delta_n$ be arbitrary. 
\begin{enumerate}
\item\label{thm2.y>} If $y>\sum_{i=1}^n\delta_i^+$, then $\bm{y}^*$ is a solution to the $\ell^1$ problem \eqref{l1.prob} if and only if 
\begin{equation*}
y_i^*= \delta_i^+ +\varepsilon_i\quad\mbox{for all}\quad  i\in[n]
\end{equation*} 
where $\varepsilon_i\ge0$ are any values summing to $y-\sum_{i=1}^n\delta_i^+$. A particular solution is $y_i^*= \delta_i^+ +\varepsilon$ for all $i\in[n]$ with $\varepsilon=(y-\sum_{i=1}^n\delta_i^+)/n$.

\item\label{thm2.y<=}  Otherwise, $y\le \sum_{i=1}^n\delta_i^+$. Then $\bm{y}^*$ is a solution to the $\ell^1$ problem \eqref{l1.prob} if and only if
\begin{equation*}
y_i^*=\alpha_i\delta_i^+\quad\mbox{for all}\quad i\in[n]
\end{equation*}
where the $\alpha_i$ are values in the hyperplane $\{\bm{\alpha}\in[0,1]^n: \bm{\alpha}'\bm{\delta}^+=y\}$.  A particular solution is $y_i^*=\alpha\delta_i^+$ for all $i\in[n]$ with $\alpha=y/\sum_{i=1}^n\delta_i^+\in(0,1]$.
\end{enumerate}
\end{theorem}

In the portfolio rebalancing problem in which the $\delta_i$ are naive adjustments~\eqref{del.def} which sum to $y$, we have $\sum_{i=1}^n\delta_i^+\ge \sum_{i=1}^n\delta_i=y$, so case~\ref{thm2.y<=} of the theorem applies. 

\begin{proof}  Let $f(\bm{y})=\sum_{i=1}^n\left|y_i-\delta_i\right|$. 

\emph{Case~\ref{thm2.y>}: $y>\sum_{i=1}^n\delta_i^+$.} We begin by showing that any $\bm{y}\in \mathbb{R}_+^n$ summing to $y$ with $y_i<\delta_i^+$ for some $i$ can be improved upon by increasing $y_i$ to at least $\delta_i^+$ while maintaining the sign of $y_j-\delta_j^+$ in the remaining components; this allows us to consider only $\bm{y}$ with $y_i\ge\delta_i^+$ for all $i$.  To this end, suppose $\bm{y}$ sums to $y$ but $y_1<\delta_1^+$. Since $y_1\ge 0$, we know that $0<\delta_1^+=\delta_1$, so let $\varepsilon:=\delta_1-y_1>0$. For $i=2,\ldots,n$ let $d_i:=(y_i-\delta_i^+)^+$. We have
\begin{multline*}
\sum_{i=2}^n d_i \ge \sum_{i=2}^n (y_i-\delta_i^+)=\left(\sum_{i=1}^n y_i - y_1\right)-\left(\sum_{i=1}^n \delta_i^+-\delta_1^+\right)\\
 = \left(y-\sum_{i=1}^n \delta_i^+ \right)+(\delta_1-y_1)  >0+\varepsilon=\varepsilon.
\end{multline*}
Let $\gamma:= \varepsilon/\sum_{i=2}^nd_i\in(0,1)$, and $\widetilde{\bm{y}}:=(\delta_1, y_2-\gamma d_2, y_3-\gamma d_3,\ldots, y_n-\gamma d_n)$. Note that $\widetilde{\bm{y}}$ sums to $y$ because
$$\sum_{i=1}^n(y_i-\widetilde{y}_i) = (y_1-\delta_1)+\gamma\sum_{i=2}^n d_i=-\varepsilon+\varepsilon=0.$$
Also note that $\mbox{sign}(\widetilde{y}_i-\delta_i^+) = \mbox{sign}(y_i-\delta_i^+)$ for all $i\ge 2$ because, if  $0<d_i=y_i-\delta_i^+$, then
\begin{equation}\label{ytild-d}
\widetilde{y}_i-\delta_i^+ = (y_i-\gamma d_i)-\delta_i^+=(1-\gamma)d_i>0
\end{equation}
as well. Otherwise, $d_i=0$ so $\widetilde{y}_i=y_i$. These same arguments also show that $\mbox{sign}(\widetilde{y}_i-\delta_i) = \mbox{sign}(y_i-\delta_i)$ for all $i\ge 2$: The $d_i>0$ case follows from \eqref{ytild-d}, and the $d_i=0$ case is again trivial.  Using these facts and that $\widetilde{y}_i<y_i$ if $d_i>0$, we finally see that  $\widetilde{\bm{y}}$ decreases $f$: 
\begin{multline*}
f(\widetilde{\bm{y}})=0+\sum_{i:\; d_i>0}(\widetilde{y}_i-\delta_i) +\sum_{i:\; d_i=0} (\delta_i-y_i)< \varepsilon+\sum_{i:\; d_i>0}(y_i-\delta_i) \\
+\sum_{i:\; d_i=0} (\delta_i-y_i)=f(\bm{y}).
\end{multline*}

With this reduction we consider only $\bm{y}$ with $y_i=\delta_i^++\varepsilon_i$ where $\varepsilon_i\ge 0$ are values that sum to $y-\sum_{i=1}^n \delta_i^+$ to ensure that $\sum_{i=1}^n y_i=y$.  Since $y_i\ge\delta_i^+\ge\delta_i$ for all $i$, for any $\bm{y}$ of this form we have
\begin{multline*}
f(\bm{y})=\sum_{i=1}^n(y_i-\delta_i)=\sum_{i:\; \delta_i>0}\varepsilon_i+ \sum_{i:\; \delta_i\le 0}(\varepsilon_i-\delta_i) = \sum_{i=1}^n\varepsilon_i- \sum_{i:\; \delta_i\le 0} \delta_i\\
 = y-\sum_{i=1}^n \delta_i^+- \sum_{i:\; \delta_i\le 0} \delta_i=y-\sum_{i=1}^n\delta_i.
\end{multline*}
This does not depend on the particular $\varepsilon_i$ so all solutions of this form minimize $f$. The particular case in the theorem is obtained by considering $\varepsilon_i$ constant.

\emph{Case~\ref{thm2.y<=}: $y\le \sum_{i=1}^n\delta_i^+$.} By similar arguments as in case~\ref{thm2.y>}, this case can be reduced to considering $\bm{y}$ with $y_i\le\delta_i^+$ for all $i$. This is done by taking a $\bm{y}$ summing to $y$ with $y_1>\delta_1^+$, and creating $\widetilde{\bm{y}}=(\delta_1^+,y_2+\gamma d_2, y_3+\gamma d_3, \ldots, y_n+\gamma d_n)$ where $\gamma$ is as above but with $\varepsilon=y_1-\delta_1^+$ and  $d_i=(\delta_i^+-y_i)^+$. By similar arguments it can be shown that $f(\widetilde{\bm{y}})<f(\bm{y})$.

With this reduction we consider only $\bm{y}$ with $y_i=\alpha_i\delta_i^+$ for all $i$, where $\alpha_i\in[0,1]$ and satisfy the defining condition of the hyperplane in the theorem, which guarantees that the $y_i$ sum to $y$. For  $\bm{y}$ of this form,
\begin{multline*}
f(\bm{y})=\sum_{i:\; \delta_i>0}(1-\alpha_i)\delta_i+ \sum_{i:\; \delta_i\le 0}|\delta_i| =\sum_{i=1}^n|\delta_i|-\sum_{i:\; \delta_i>0} \alpha_i\delta_i \\
=\sum_{i=1}^n|\delta_i|-\sum_{i=1}^n \alpha_i\delta_i^+ =\sum_{i=1}^n|\delta_i|-y,
\end{multline*}
which does not depend on the particular $\alpha_i$, so  any $\bm{y}$ of this form achieves the minimum.  The particular case in the theorem is found by taking $\alpha_i$ constant.
\end{proof}

\section{Examples}

\subsection{Portfolio rebalancing} 
\begin{figure}[!ht]
    \centering
    \caption{Allocations, amounts, and adjustments of a hypothetical $\$10,000$ portfolio, rounded to the nearest \$1 and 1\%.} \label{table:port}
\includegraphics[width=13cm]{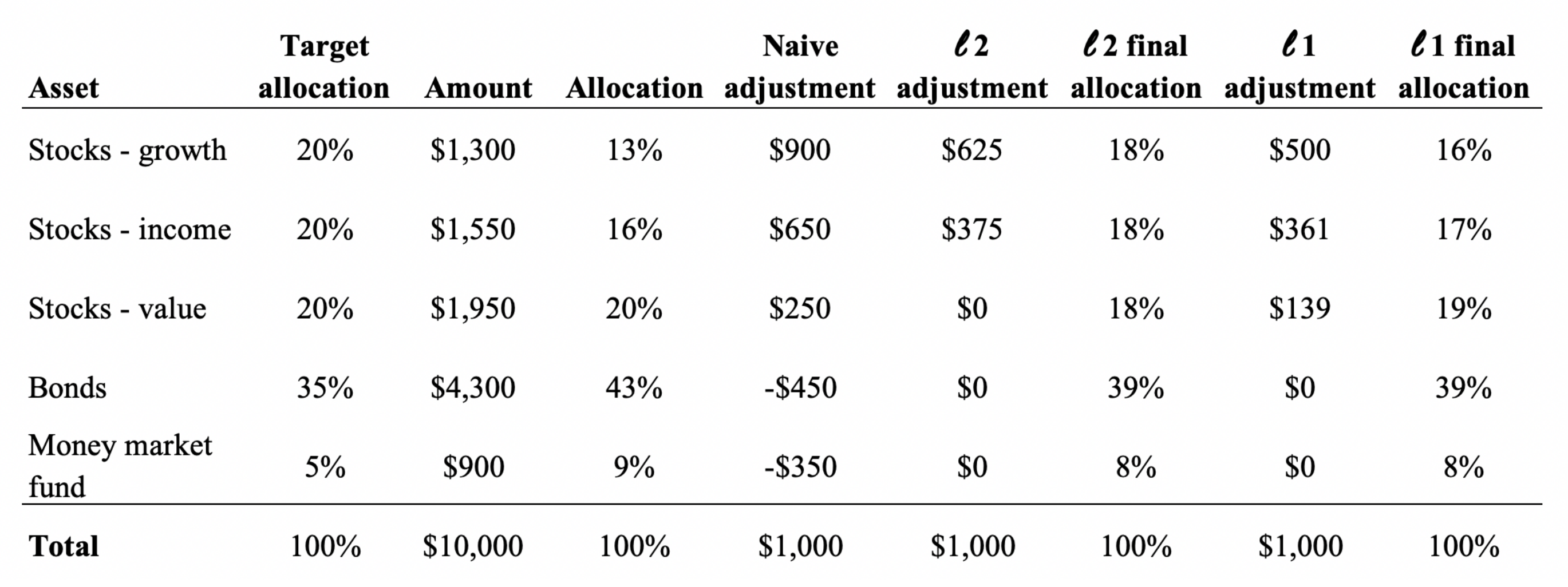} 
\end{figure}

How would you add \$1,000 to the \$10,000 portfolio in Table~\ref{table:port} to bring the assets as close as possible to the target allocation?  The hypothetical portfolio is composed of the $n=5$ assets  -- so-called growth, income, and value stocks, bonds, and a money market fund --  in the dollar amounts~$x_i$ given in the table, whose allocation differs from the target allocation~$p_i$ in the table's 2nd column, say, because of recent declines in the stock market. The naive adjustments, given by \eqref{del.def}, would require transactions on all 5 assets, including  liquidating some of the bonds and  money market fund.

\renewcommand{\figurename}{Table}

To compute the $\ell^2$ adjustments we note that, in the notation of Theorem~\ref{thm:l2}, $k^*=2$ since $$\sum_{i=1}^2 (\delta_i-\delta_2)=\$350 <\underbrace{\$1000}_{=y}<\$1050=\sum_{i=1}^3 (\delta_i-\delta_3).$$ Thus only the $k^*=2$ assets  with the largest naive adjustments -- growth and income stocks -- will be added to in the amounts $(\delta_i-\lambda^*)^+$, where the threshold is $$\lambda^*=\frac{\sum_{i=1}^{k^*}\delta_i-y}{k^*}=\frac{\$900+\$650-\$1000}{2} = \$275.$$ The positive parts over this threshold are the $\ell^2$ adjustments.

For the $\ell^1$ adjustments, the particular solution given in Part~\ref{thm2.y<=} of Theorem~\ref{thm:l1} is $\alpha\delta_i^+$ where $$\alpha=\frac{y}{\sum_{i=1}^n\delta_i^+} = \frac{\$1000}{\$900+\$650+\$250} = \frac{5}{9}=.555\ldots$$ is the deflation factor applied to all positive naive adjustments. Thus the $\ell^1$ adjustments add to all 3 stock types.

The $\ell^2$ and $\ell^1$ final allocations $(x_i+y_i^*)/(x+y)$ (with $x=$\$10,000) are what would result from making those respective adjustments, and show the difference in these approaches.  The $\ell^2$ norm penalizes larger deviations from the target more than the $\ell^1$ adjustments, and thus moves  more funds  and causes the largest shift in final allocation to the asset with  the largest deviation -- growth stocks --  and does not add to the value stocks, which have the smallest positive naive adjustment.   On the other hand, the $\ell^1$ norm penalizes the 3 stock types' deviations more evenly and adds more to the value stocks and less to the growth stocks.  The result is a final allocation of value stocks closer to the target than the $\ell^2$ approach, but still skewed toward value stocks like the initial allocation.

\subsection{A statistical application: Maximum likelihood estimation on a simplex}\label{sec:lava} Hot magma spews from inside the Earth and lava explodes violently into the atmosphere, falling back to the Earth's surface and eventually cooling into basalt rock. Geologists gain insight into this process by studying the chemical makeup of the basalt. \cite{Thompson72} gives an example of one such study on the Isle of Skye in Scotland where basalt samples were chemically analyzed resulting in data vectors like

\scriptsize
\begin{center}
\begin{tabular}{ccccccccccc}
$\mathrm{SiO_2}$ & $\mathrm{Al_2O_3}$ &$\mathrm{Fe_2O_3}$ &$\mathrm{MgO}$ &$\mathrm{CaO}$ &$\mathrm{Na_2O}$ &$\mathrm{K_2O}$ &$\mathrm{TiO_2}$ &$\mathrm{P_2O_5}$ &$\mathrm{MnO}$ & Other\\
46.31\%&14.18\%&12.32\% &12.74\% &9.62\% &2.51\% &0.34\% &1.53\% &0.16\% &0.18\% &0.11\% 
\end{tabular}
\end{center}
\normalsize

From vectors like this, the task at hand is to estimate the true percentages of these compounds in the Isle of Skye lava. Estimating a vector of percentages, which sum to $1$ and hence live on a simplex, fall in the area of statistics called \textit{compositional data analysis} \cite{Aitchison82} and arise in many areas beyond geology including economics, demographics, and medicine. 

Continuing with the lava example, volcanologists may draw samples from various locations near a volcanic eruption resulting in multiple vectors like that one, each with varying percentages of the compounds due to natural variation, measurement error, and other factors. Let $\theta_1,\ldots,\theta_{11}$ denote the true overall percentages of the 11 categories above (10 chemical compounds, plus ``Other'') and suppose $z_1,\ldots,z_{11}$ are the averages of many data vectors like that one. A reasonable statistical model for the natural chemical variation in the basalt as well as the inherent measurement error in the chemical analyses is to suppose that the $z_i$ are independent Gaussian random variables with the correct mean $\theta_i$ and some inherent variance~$\sigma^2$,
$$z_i\sim N(\theta_i,\sigma^2).$$
After writing out the likelihood function of the $z_i$ and taking logs, the maximum likelihood estimate of the $\theta_i$ is then the solution of the following problem:
\begin{equation}\label{mle}
\mbox{find $\theta_1,\ldots,\theta_{11}\ge 0$ minimizing}\quad \sum_{i=1}^{11}\left(\theta_i -z_i \right)^2 \quad\mbox{subject to}\quad  \sum_{i=1}^{11} \theta_i=1.
\end{equation} 
Although this problem may not appear to have anything in common with the resource allocation problems mentioned above, \eqref{mle} is obviously a special case of \eqref{l2.prob} with $\theta_i$ playing the role of $y_i$, $z_i$ playing the role of $\delta_i$, and $y=1$. Theorem~\ref{thm:l2} can be  applied directly. 

\begin{abstract}
How do you bring your assets as close as possible to your target allocation by only investing a fixed amount of additional funds, and not selling any assets?  We look at two versions of this problem which have simple, closed form solutions revealed by basic calculus and algebra.
\end{abstract}

\vfill\eject

\end{document}